\documentstyle[twoside]{article}
\title{On a new result for the hypergeometric function}
\author{\sc Arjun K. Rathie$^a$ and Richard B. Paris$^b$\\
\\
${}^a\!$ Department of Mathematics, Vedant College of Engineering and Technology \\
(Rajasthan Technical University), Bundi, 323021, Rajasthan, India\\
E-Mail: arjunkumarrathie@gmail.com\\
${}^b\!$ Division of Computing and Mathematics,\\
Abertay University, Dundee DD1 1HG, UK\\
E-Mail: r.paris@abertay.ac.uk}
\begin{document}
\def\f#1#2{\mbox{${\textstyle \frac{#1}{#2}}$}}
\def\dfrac#1#2{\displaystyle{\frac{#1}{#2}}}
\def\boldal{\mbox{\boldmath $\alpha$}}
\newcommand{\bee}{\begin{equation}}
\newcommand{\ee}{\end{equation}}
\newcommand{\lam}{\lambda}
\newcommand{\ka}{\kappa}
\newcommand{\al}{\alpha}
\newcommand{\th}{\theta}
\newcommand{\om}{\omega}
\newcommand{\Om}{\Omega}
\newcommand{\fr}{\frac{1}{2}}
\newcommand{\fs}{\f{1}{2}}
\newcommand{\g}{\Gamma}
\newcommand{\br}{\biggr}
\newcommand{\bl}{\biggl}
\newcommand{\ra}{\rightarrow}
\renewcommand{\topfraction}{0.9}
\renewcommand{\bottomfraction}{0.9}
\renewcommand{\textfraction}{0.05}
\newcommand{\mcol}{\multicolumn}
\date{}
\maketitle
\pagestyle{myheadings}
\markboth{\hfill  {\it A.K. Rathie and R.B. Paris }  \hfill}
{\hfill {\it A hypergeometric function relation} \hfill}
\begin{abstract}
The aim of this note is to provide a new identity connected with the Gauss hypergeometric function.
This is achieved using results of certain combinatorial identities and a hypergeometric function approach.
\vspace{0.4cm}

\setcounter{equation}{0}
\vspace{0.3cm}

\end{abstract}
\vspace{0.3cm}

\begin{center}
{\bf 1. \  Introduction}
\end{center}
\setcounter{section}{1}
\setcounter{equation}{0}
\renewcommand{\theequation}{\arabic{section}.\arabic{equation}}
In probability theory and combinatorial identities, the so-called match box problem is well known \cite{WF}. This may be stated as follows:
\begin{quote}
If matches are drawn one at a time and at random from two match boxes, each initially containing $n$ matches, then the probability $P_{n,r}$ that when one box is found empty, the other will contain exactly $r$ matches is given by
\[P_{n,r}=2^{r-2n}\bl(\!\!\!\begin{array}{c}2n-r\\n\end{array}\!\!\!\br),\quad r=0, 1, 2, \ldots , n.\]
\end{quote}
Feller \cite{WF} obtained the mean and Riordan \cite{Ri} obtained the variance for this problem by rather lengthy methods. One of the present authors \cite{AKR1} found general moments and in particular the mean and variance by a simple combinatorial method.

The above match box problem was generalised by Rohatgi \cite{Ro} as follows:
\begin{quote}
Let matches be drawn one at a time and at random from two match boxes, the first box is selected with probability $p$ and the second box with probability $q=1-p$, each box initially containing $n$ matches. Then the probability $P_r(n,p)$ that when one box is found empty, the other will contain exactly $r$ matches is given by
\bee\label{e12}
P_r(n,p)=\bl(\!\!\!\begin{array}{c}2n-r\\n\end{array}\!\!\!\br)p^{n+1}q^{n-r}+\bl(\!\!\!\begin{array}{c}2n-r\\n\end{array}\!\!\!\br)q^{n+1} p^{n-r}
\ee
for $r=0, 1, 2, \ldots , n.$
\end{quote}

In 2005, Rathie and Rathie \cite{AKR2} obtained expressions for the general moments, the moment generating function and the probability generating function for this problem and derived their results in terms of the readily computable Gauss hypergeometric function ${}_2F_1(a,b;c,x)$ \cite[p.~384]{DLMF}. In this note we provide a (possibly)) new result expressing the sum of two such hypergeometric functions of argument $p^{-1}$ and $q^{-1}$ by evaluating $\sum_{r=0}^nP_r(n,p)$ in two ways. This is achieved with the help of results of certain combinatorial identities and a hypergeometric function approach.

\vspace{0.6cm}

\begin{center}
{\bf 2. \  Main result}
\end{center}
\setcounter{section}{2}
\setcounter{equation}{0}
\renewcommand{\theequation}{\arabic{section}.\arabic{equation}}
The new hypergeometric function result to be established is given in the following theorem:
\newtheorem{theorem}{Theorem}
\begin{theorem}
For non-negative integer $n$ and $0<p, q<1$ with $p+q=1$, the following result holds true:
\bee\label{e21}
p\,{}_2F_1\bl(\!\!\begin{array}{c}-n,\ 1\\-2n\end{array}\!\!; \frac{1}{q}\br)+q\,{}_2F_1\bl(\!\!\begin{array}{c}-n,\ 1\\-2n\end{array}\!\!; \frac{1}{p}\br)=\frac{(n!)^2}{p^nq^n (2n)!}=\frac{(1)_n}{2^{2n}p^nq^n (\fs)_n},
\ee
where $(a)_n=\g(a+n)/\g(a)$ is the Pochhammer symbol.
\end{theorem}

\noindent {\it Proof.}\ \ Since $P_r(n,p)$ defined in (\ref{e12}) denotes a probability, then it follows that the sum of the probabilities over $r\leq n$ must be unity; that is $\sum_{r=0}^n P_r(n,p)=1$.

We now demonstrate this fact algebraically. Let
\bee\label{e22}
S_n=\sum_{r=0}^nP_r(n,p)=S_n^{(1)}+S_n^{(2)},
\ee
where, from (\ref{e12}),
\bee\label{e20}
S_n^{(1)}=\sum_{r=0}^n \bl(\!\!\!\begin{array}{c}2n-r\\n\end{array}\!\!\!\br) p^{n+1}q^{n-r},\qquad
S_n^{(2)}=\sum_{r=0}^n \bl(\!\!\!\begin{array}{c}2n-r\\n\end{array}\!\!\!\br) q^{n+1}p^{n-r}.
\ee

Making the change of summation index $r\to n-r$, we easily see that
\begin{eqnarray*}
S_n^{(1)}&=&\sum_{r=0}^n \bl(\!\!\!\begin{array}{c}n+r\\n\end{array}\!\!\!\br) p^{n+1}q^{r}
=\sum_{r=0}^n\bl\{\bl(\!\!\!\begin{array}{c}n+r-1\\n\end{array}\!\!\!\br)+\bl(\!\!\!\begin{array}{c}n+r-1\\n-1\end{array}\!\!\!\br)\br\} p^{n+1}q^r\\
&=&q\bl\{S_n^{(1)}-\bl(\!\!\!\begin{array}{c}2n\\n\end{array}\!\!\!\br) p^{n+1}q^n\br\}+p\bl\{S_{n-1}^{(1)}+\bl(\!\!\!\begin{array}{c}2n-1\\n-1\end{array}\!\!\!\br)p^nq^n\br\}.
\end{eqnarray*}
From this it follows that, since $p+q=1$,
\[S_n^{(1)}=S_{n-1}^{(1)}+p^nq^n \bl(\!\!\!\begin{array}{c}2n\\n\end{array}\!\!\!\br) \bl(\frac{1}{2}-q\br).\]
In a similar manner we obtain
\[S_n^{(2)}=S_{n-1}^{(2)}+p^nq^n \bl(\!\!\!\begin{array}{c}2n\\n\end{array}\!\!\!\br) \bl(\frac{1}{2}-p\br).\]

Thus, it easily seen that
\[S_n=S_{n-1}= \ldots =S_0=1\]
and hence that
\bee\label{e23}
S_n=\sum_{r=0}^nP_r(n,p)=1.
\ee

Again, from (\ref{e21}) and (\ref{e20}), we have
\[S_n=\sum_{r=0}^n\bl(\!\!\!\begin{array}{c}2n-r\\n\end{array}\!\!\!\br)p^{n+1}q^{n-r}+\sum_{r=0}^n\bl(\!\!\!\begin{array}{c}2n-r\\n\end{array}\!\!\!\br)q^{n+1} p^{n-r}.\]
Using the results
\[\bl(\!\!\!\begin{array}{c}n\\r\end{array}\!\!\!\br)=\frac{\g(n+1)}{\g(r+1) \g(n+1-r)},\qquad \g(\alpha-r)=\frac{(-1)^r \g(\alpha)}{(1-\alpha)_r},\]
we have after some algebra,
\[S_n=p^nq^n \frac{(2n)!}{(n!)^2}\,\bl\{p\sum_{r=0}^n\frac{(-n)_r(1)_r}{(-2n)_r r!}\,q^{-r}+q\sum_{r=0}^n\frac{(-n)_r(1)_r}{(-2n)_r r!}\,p^{-r}\br\}.\]
The sums can be expressed as terminating Gauss hypergeometric functions to yield
\bee\label{e24}
S_n=p^nq^n \frac{(2n)!}{(n!)^2}\,\bl\{p\,{}_2F_1\bl(\!\!\begin{array}{c}-n,\ 1\\-2n\end{array}\!\!; \frac{1}{q}\br)+q\,{}_2F_1\bl(\!\!\begin{array}{c}-n,\ 1\\-2n\end{array}\!\!; \frac{1}{p}\br)\br\}.
\ee

Then on equating (\ref{e23}) and (\ref{e24}), we arrive at the desired result (\ref{e21}) asserted in the theorem. This completes the proof of (\ref{e21}).   $\Box$
\medskip

\newtheorem{corollary}{Corollary}
\begin{corollary} \ \ In (\ref{e21}), if we set $p=q=\fs$, we find
\[{}_2F_1\bl(\!\!\begin{array}{c}-n, 1\\-2n\end{array}\!\!;1\br)=\frac{(1)_n}{(\fs)_n},\]
which is a well-known result in the literature of hypergeometric functions recorded, for example, in \cite{PBM}.
\end{corollary}
\medskip

\noindent{\bf Remark.}\ \ The evaluation of ${}_2F_1(a,b;c;x)$ when $c$ and $a$ (or $b$) are negative integers needs some care. When $c$ and $a$ (or $b$) are independent of each other, as in (\ref{e21}) where $c=-2n$, $a=-n$ and $b=1$, the series expansion of this function terminates after $n+1$ terms. If, however, $c$ and $a$ (or $b$) are connected then the function consists of a finite sum with $n+1$ terms together with an infinite series; see \cite[p.~109]{T}.

\end{document}